\documentclass[11pt,reqno]{amsart}
\usepackage[left=33mm,right=33mm,top=30mm,bottom=32mm]{geometry}
\usepackage{mathtools,amssymb,amsthm,mathrsfs,color,float}
\usepackage{paralist}
\usepackage{stackengine}

\usepackage[colorlinks,
linkcolor=red,
anchorcolor=green,
citecolor=blue, 
]{hyperref}

\usepackage[T1]{fontenc}
\usepackage[utf8]{inputenc}

\usepackage{calc}
\definecolor{bleu1}{RGB}{0,57,128}
\def\bleu1{\color{bleu1}}

\usepackage{etoolbox}
\patchcmd{\section}{\normalfont}{\normalfont \bleu1}{}{}
\patchcmd{\subsection}{\normalfont}{\normalfont \bleu1}{}{}
\patchcmd{\subsubsection}{\normalfont}{\normalfont \bleu1}{}{}

\newtheorem{The}{\bleu1 Theorem}[section]
\newtheorem{Pro}{\bleu1 Proposition}[section]
\newtheorem{Lem}{\bleu1 Lemma}[section]

\theoremstyle{definition}

\newcommand{\T}{\mathbb{T}}
\newcommand{\R}{\mathbb{R}}

\newcommand{\N}{\mathbb{N}}

\newcommand\Set[2]{ 
\left\{#1\ \middle\vert\ #2 \right\}}

\makeatletter
\@namedef{subjclassname@2020}{\textup{2020} Mathematics Subject Classification}
\makeatother

\title[Qualitative Estimates of Topological Entropy]{Qualitative Estimates of Topological Entropy for Non-Monotone Contact Lax-Oleinik Semiflow}
\author{Wei Cheng, Jiahui Hong, Zhi-Xiang Zhu$^{\ast}$}
\address[Wei Cheng]{School of Mathematics, Nanjing University, Nanjing 210093, China}
\email{chengwei@nju.edu.cn}
\address[Jiahui Hong]{School of Mathematics, Nanjing University of Aeronautics and Astronautics, Nanjing 211106, China}
\email{hongjiahui@nuaa.edu.cn}
\address[Zhi-Xiang Zhu]{School of Mathematics, Nanjing University, Nanjing 210093, China}
\email{zhuzhixiang0928@gmail.com}
\thanks{$^{\ast}$ Corresponding author.}

\begin{document}
\maketitle

\begin{abstract}
For the non-monotone Hamilton-Jacobi equations of contact type, the associated Lax-Oleinik semiflow \((T_t, C(M))\) is expansive. In this paper, we provide qualitative estimates for both the lower and upper bounds of the topological entropy of the semiflow.
\end{abstract}


\section{Introduction}

Let $M$ be a compact and connected manifold without boundary. A function $H=H(x,p,u):T^* M\times\R\to\R$ is called a \emph{Tonelli Hamiltonian of contact type} if it is of class $C^2$ and satisfies the following conditions:
\begin{enumerate}[\rm (H1)]
	\item $H_{pp}(x,p,u)>0$, $\forall (x,p,u)\in T^*M\times\R$;
	\item There exist a super-linear continuous function $\theta^*_1:[0,+\infty)\to[0,+\infty)$ and a constant $c_0\geqslant0$ such that
	\begin{align*}
		H(x,p,0)\geqslant\theta^*_1(|p|_x)-c_0,\qquad\forall(x,p)\in T^*M;
	\end{align*}
	\item There exists $\kappa>0$ such that
	\begin{align*}
		|H_u(x,p,u)|\leqslant \kappa,\qquad \forall (x,p,u)\in T^*M\times\R.
	\end{align*}
\end{enumerate}
Consider the Hamilton-Jacobi equation of contact type:
\begin{equation}\label{eq:HJe_intr}
	\begin{cases}
		D_tu+H(x,Du(x),u(x))=0, &\quad t>0, x\in M,\\
		u(0,x)=\phi(x), &\quad x\in M, 
	\end{cases}
\end{equation}
where \( \phi\in C(M) \), the space of continuous functions on \( M \). Using Herglotz' variational principle (\cite{CCWY2019,CCJWY2020}) or the implicit variational principle (\cite{Wang_Wang_Yan2017}), the (unique) solution of \eqref{eq:HJe_intr} has a representation formula. More precisely,
\begin{align*}
	u(t,x) = T_t\phi(x), \quad (x \in M, t \geqslant 0).
\end{align*}

The Lax-Oleinik operators \( \{T_t\} \) are defined as follows:
\begin{align*}
	(T_t\phi)(x) = \inf_{\xi \in \mathcal{A}_{t,x}} \left\{ \phi(\xi(0)) + \int_0^t L(\xi(s), \dot{\xi}(s), u_{\xi}(s)) \, ds \right\},
\end{align*}
where $L:TM\times\R\to\R$ is the associated Tonelli Lagrangian of contact type defined by
\begin{align*}
	L(x,v,u):=\sup_{p\in T^*_xM}\{p(v)-H(x,p,u)\},\qquad x\in M, v\in T_xM, u\in\R,
\end{align*}
$\mathcal{A}_{t,x}$ is the family of absolutely continuous arcs $\xi:[0,t]\to M$ with $\xi(t)=x$, and \( u_{\xi} : [0,t] \to \mathbb{R} \) uniquely solves the Carath\'eodory equation:
\begin{align*}
	\begin{cases}
		\dot{u}_{\xi}(s) = L(\xi(s), \dot{\xi}(s), u_{\xi}(s)), & \quad s \in [0,t], \\
		u_{\xi}(0) = \phi(\xi(0)). &
	\end{cases}
\end{align*}
For more details on the semigroup \( \{T_t\}_{t \geqslant 0} \), see Section \ref{sec:Herglotz}. A typical example of a contact-type Hamilton-Jacobi equation is the discounted equation with \( H(x,p,u) = \lambda u + H(x,p) \), where \( \lambda > 0 \).

To understand the dynamics of the semiflow $(T_t,C(M))$ generated by \eqref{eq:HJe_intr}, especially when \( \lambda < 0 \) in the discounted equation or \( H_u < 0 \) in general, many authors have studied the dynamics in a deterministic way, such as through Lyapunov stability and periodicity (\cite{Jin_Yan2021,Wang_Yan_Zhao2023}). In the current paper, we will concentrate on certain statistical behaviors, specifically the topological entropy $h_{\rm top}(T_t)$ (see Section \ref{sec:entropy} for the definition) of the semiflow \( \{T_t\} \). Recall that one can easily deduce $h_{\rm top}(T_t)=0$ for the classical Lax-Oleinik semigroup from the qualitative estimate of the Kolmogorov \( \varepsilon \)-entropy in the space \( W^{1,1}(M) \), which was obtained in \cite{Ancona_Cannarsa_Nguyen2016_1, Ancona_Cannarsa_Nguyen2016_2}. 

When discussing the contact-type Lax-Oleinik semigroup, we contrast the case \( h_{\rm top}(T_t) = 0 \) when \( H_u > \lambda > 0 \), due to the contraction property of \( \{T_t\} \), with the expansion case when \( H_u < 0 \). In the latter case, we will primarily focus on giving a qualitative estimate of \( h_{\rm top}(T_t) \). 

For the lower bound of $h_{\rm top}(T_t)$, we have the following result: under the assumptions (H1)-(H3), if there exists $0 < \lambda \leqslant \kappa$ such that $\lambda \leqslant L_u(x,v,u) \leqslant \kappa$ for all $(x,v,u) \in TM \times \mathbb{R}$, then $h_{\rm top}(T_t) \geqslant \lambda$ (see Theorem \ref{thm:entrop2}). The estimate of the upper bound for $h_{\rm top}(T_t)$ is more difficult. We will focus on the discounted Hamiltonian in the form
\begin{align*}
	H(x,p,u) = h(p) - \lambda u, \quad (x,p,u) \in T^*\T^d \times \mathbb{R},
\end{align*}
where $h = h(p): T^*\T^d \to \mathbb{R}$ is a Tonelli Hamiltonian and $\lambda > 0$ is a constant. We prove that for such a Hamiltonian, $h_{\rm top}(T_t) = \lambda$ (Theorem \ref{thm:dis}).

We remark that it is still an interesting problem whether the developed theory in infinite-dimensional dynamical systems (see, for instance, \cite{Lian_Liu_Lu2016, Lu_Wang_Young2013, Lian_Lu2010, Huang_Lu2017}) can be applied to study the expansive Lax-Oleinik semiflow \((T_t, C(M))\) for \eqref{eq:HJe_intr}, with positive topological entropy, and its potential consequences on the complicated dynamical behavior.

\medskip

\noindent\textbf{Acknowledgements.} 
Wei Cheng was partly supported by National Natural Science Foundation of China (Grant No. 12231010). Wei Cheng and Jiahui Hong are grateful to Kaizhi Wang for helpful discussions.

\section{Preliminaries}

\subsection{Topological entropy}\label{sec:entropy}

Let \( (X,d) \) be a complete metric space, and let \( \{T_t\}_{t \geqslant 0} \) be a continuous semigroup action on \( X \). That is, \( T(t,x) = T_t(x): [0,+\infty) \times X \to X \) is continuous, \( T_0 = \text{id} \), and \( T_s \circ T_t = T_{s+t} \) for any \( t, s \geqslant 0 \). For any \( x \in X \), \( t \geqslant 0 \), and \( \varepsilon > 0 \), we define the \emph{dynamical ball} of center \( x \) with length \( t \) and radius \( \varepsilon > 0 \) as
\begin{align*}
	B(x,t,\varepsilon) = \Set{y \in X}{d(T_s(y), T_s(x)) < \varepsilon \ \text{for all} \ s \in [0,t]}.
\end{align*}
It is immediate that \( B(x,0,\varepsilon) = B(x,\varepsilon) \). Let \( K \subset X \) be a compact subset. A set \( E \subset X \) is called a \( (t,\varepsilon) \)-generating set for \( K \) if
\begin{align*}
	K \subset \bigcup_{x \in E} B(x,t,\varepsilon),
\end{align*}
and \( E \subset K \) is a \( (t,\varepsilon) \)-separated set if the dynamical ball \( B(x,t,\varepsilon) \) for each \( x \in E \) contains no other element of \( E \). For any compact subset \( K \) and \( \varepsilon > 0 \), let \( g_t(\varepsilon,K) \) be the smallest cardinality of a \( (t,\varepsilon) \)-generating set of \( K \), and \( s_t(\varepsilon,K) \) the largest cardinality of a \( (t,\varepsilon) \)-separated set \( E \subset K \). Then, we define
\begin{align*}
	g(K) &= \lim_{\varepsilon \to 0^+} g(K,\varepsilon) = \lim_{\varepsilon \to 0^+} \limsup_{t \to \infty} \frac{1}{t} \log g_t(\varepsilon,K), \\
	s(K) &= \lim_{\varepsilon \to 0^+} s(K,\varepsilon) = \lim_{\varepsilon \to 0^+} \limsup_{t \to \infty} \frac{1}{t} \log s_t(\varepsilon,K).
\end{align*}
It is known that \( g(K) = s(K) \) for any compact set \( K \subset X \). We now define the topological entropy of the semigroup \( \{T_t\} \) as
\begin{align*}
	h(T_t) = h_{\rm top}(T_t) = \sup_{K \subset X} g(K) = \sup_{K \subset X} s(K),
\end{align*}
where the supremum is taken over all compact subsets \( K \) of \( X \).

Recall a semigroup $\{T_t\}_{t\geqslant0}$ is called \emph{uniformly continuous} if for any $t>0$ and $\varepsilon>0$ there exists $\delta>0$ such that
\begin{align*}
	d(x,y)<\delta\Longrightarrow d(T_s(x),T_s(y))<\varepsilon\ \text{for every}\ s\in[0,t].
\end{align*}

\begin{Pro}[\cite{Viana_Oliveira_book2016}]\label{pro:discrete}
If the semi-flow $\{T_t\}$ is uniformly continuous then its topological entropy $h_{\rm top}(T_t)$ coincides with the topological entropy of its time-$\tau$ map for any $\tau>0$.	
\end{Pro}

\begin{Lem}[\cite{Viana_Oliveira_book2016}]\label{lem:entropy1}
If the semigroup $\{T_t\}$ is non-expansive, i.e.,
	\begin{align*}
		d(T_t(x),T_t(y))\leqslant d(x,y),\qquad\forall t\geqslant0, x,y\in X, 
	\end{align*}
then for any compact subset $K$ of $X$, $g(K)=0$. In particular, $h_{\rm top}(T_t)=0$.
\end{Lem}

In fact, the topological entropy \( h_{\rm top}(T_t) \) is independent of the action of the semigroup \( T_t \) in finite time. More precisely, for any \( x \in X \), \( t \geqslant t_0 > 0 \), and \( \varepsilon > 0 \), we define the dynamical ball over time \( [t_0,t] \) as
\begin{align*}
	B_{t_0}(x,t,\varepsilon) = \Set{y \in X}{d(T_s(x), T_s(y)) < \varepsilon \text{ for all } s \in [t_0,t]}.
\end{align*}
For any compact subset \( K \subset X \), \( t \geqslant t_0 > 0 \), and \( \varepsilon > 0 \), let
\begin{align*}
	g_{t_0,t}(\varepsilon,K) = \min\Set{ \#(E)}{E \subset X, \ \bigcup_{x \in E} B_{t_0}(x,t,\varepsilon) \supset K}.
\end{align*}

\begin{Pro}\label{pro:cut}
\hfill
\begin{enumerate}[\rm (1)]
	\item For any compact subset $K\subset X$, $t\geqslant t_0>0$, $\varepsilon>0$, we have
		\begin{align*}
			g_{t_0,t}(\varepsilon,K)\leqslant g_{t}(\varepsilon,K),\quad g_{t_0}(\varepsilon,K)\cdot g_{t_0,t}(\varepsilon,K) \geqslant g_{t}(2\varepsilon,K).
		\end{align*}
	\item For any compact subset $K\subset X$, $t_0>0$, we have
		\begin{align*}
			g(K)=\lim_{\varepsilon\to 0^+}\limsup_{t\to+\infty}\frac{1}{t}\log g_{t_0,t}(\varepsilon,K).
		\end{align*}
\end{enumerate}
\end{Pro}

\begin{proof}
By definition, we know that
\begin{align*}
	B(x,t,\varepsilon) \subset B_{t_0}(x,t,\varepsilon), \quad \forall x \in X,
\end{align*}
which implies directly
\begin{align*}
    g_{t_0,t}(\varepsilon,K) \leqslant g_{t}(\varepsilon,K).
\end{align*}
For the second inequality, we consider two finite covers of the compact set \( K \):
\begin{align*}
	\bigcup_{i=1}^{g_{t_0}(\varepsilon,K)} B(x_i,t_0,\varepsilon) \supset K, \quad 
	\bigcup_{j=1}^{g_{t_0,t}(\varepsilon,K)} B_{t_0}(y_j,t,\varepsilon) \supset K.
\end{align*}
If \( B(x_i,t_0,\varepsilon) \cap B_{t_0}(y_j,t,\varepsilon) \neq \varnothing \), we choose \( z_{ij} \in B(x_i,t_0,\varepsilon) \cap B_{t_0}(y_j,t,\varepsilon) \). It follows from the definition of the dynamical ball that
\begin{align*}
	B(z_{ij},t,2\varepsilon) \supset B(x_i,t_0,\varepsilon) \cap B_{t_0}(y_j,t,\varepsilon).
\end{align*}
Thus, we have
\begin{align*}
	g_{t}(2\varepsilon,K) \leqslant \#\{z_{ij}\} \leqslant g_{t_0}(\varepsilon,K) \cdot g_{t_0,t}(\varepsilon,K).
\end{align*}

Now we turn to prove statement (2). Using the result in (1), we obtain
\begin{align*}
	g(K) &= \lim_{\varepsilon \to 0^+} \limsup_{t \to +\infty} \frac{1}{t} \log g_{t}(2\varepsilon,K) \\
	&\leqslant \lim_{\varepsilon \to 0^+} \limsup_{t \to +\infty} \frac{1}{t} \log \left(g_{t_0}(\varepsilon,K) \cdot g_{t_0,t}(\varepsilon,K)\right) \\
	&= \lim_{\varepsilon \to 0^+} \limsup_{t \to +\infty} \frac{1}{t} \log g_{t_0,t}(\varepsilon,K) \\
	&\leqslant \lim_{\varepsilon \to 0^+} \limsup_{t \to +\infty} \frac{1}{t} \log g_{t}(\varepsilon,K) = g(K).
\end{align*}
This completes the proof.
\end{proof}

\subsection{Herglotz' variational principle}\label{sec:Herglotz}

Let $M$ be a compact and connected manifold without boundary and $H:T^*M\times\R\to\R$ a Tonelli Hamiltonian of contact type. Its associated Tonelli Lagrangian of contact type $L=L(x,v,u):TM\times\R\to\R$ is of class $C^2$ and satisfies the following conditions:
\begin{enumerate}[\rm (L1)]
	\item $L_{vv}(x,v,u)>0$, $\forall (x,v,u)\in TM\times\R$;
	\item There exist a super-linear continuous function $\theta_0:[0,+\infty)\to[0,+\infty)$ and a constant $c_0\geqslant0$ such that
	\begin{align*}
		L(x,v,0)\geqslant\theta_0(|v|_x)-c_0,\qquad\forall(x,v)\in TM;
	\end{align*}
	\item There exists $\kappa>0$ such that
	\begin{align*}
		|L_u(x,v,u)|\leqslant \kappa,\qquad \forall (x,v,u)\in TM\times\R.
	\end{align*}
\end{enumerate}

To introduce the Lax-Oleinik evolution of contact type, let us recall the relevant Herglotz' variational principle. For any $x,y\in M$ and $t>0$, let $\Gamma^t_{x,y}$ be the set of absolutely continuous arcs $\xi:[0,t]\to M$ connecting $x$ to $y$, and $\mathcal{A}_{t,x}$ the set of absolutely continuous arcs $\xi:[0,t]\to M$ such that $\xi(t)=x$. We define the Lax-Oleinik evolution of contact type as a operator $T_t$, $t>0$, acting on a function $\phi:M\to\R$, as
\begin{align*}
	(T_t\phi)(x)=\inf_{\xi\in\mathcal{A}_{t,x}}\Big\{\phi(\xi(0))+\int^t_0L(\xi(s),\dot{\xi}(s),u_{\xi}(s))\ ds\Big\}, 
\end{align*}
where $u_{\xi}:[0,t]\to\R$ uniquely solves the Carath\'eodory equation
\begin{equation}\label{eq:cara}
	\begin{cases}
		\dot{u}_{\xi}(s)=L(\xi(s),\dot{\xi}(s),u_{\xi}(s)),\qquad s\in[0,t],\\
		u_{\xi}(0)=\phi(\xi(0)).
	\end{cases}
\end{equation}

Let $C(M)$ the space of continuous real-valued functions on $M$ endowed with the norm
\begin{align*}
	\|f\|=\max_{x\in M}|f(x)|.
\end{align*}
The space $(C(M),\|\cdot\|)$ is a complete metric space. We can conclude from the following Lemma \ref{lem:LO_semigroup} that
\begin{align*}
	T_t: C(M)\to C(M),\qquad \phi\mapsto T_t\phi
\end{align*}
is a continuous semigroup with $T_0=\text{id}$. For the proof, one can refer to \cite{CCWY2019,CCJWY2020}.
\begin{Lem}\label{lem:LO_semigroup}
\hfill
\begin{enumerate}[\rm (1)]
	\item For any $\phi\in C(M)$, $t>0$ and $x\in M$ there exists $\xi\in\mathcal{A}_{t,x}$ such that
	\begin{align*}
		T_t\phi(x)=\phi(\xi(0))+\int^t_0L(\xi(s),\dot{\xi}(s),u_{\xi}(s))\ ds=u_{\xi}(t),
	\end{align*}
	where $u_{\xi}$ is uniquely determined by \eqref{eq:cara}.
	\item For any $\phi\in C(M)$, $(t,x)\mapsto T_t\phi(x)$ is a continuous real-valued function on $[0,+\infty)\times M$, and it is locally semiconcave on $(0,+\infty)\times M$.
	\item $T_s\circ T_t\phi=T_{s+t}\phi$ for all $s,t\geqslant0$ and $\phi\in C(M)$.
	\item We have
	\begin{align*}
		\|T_t\phi-T_t\psi\|\leqslant e^{\kappa t}\|\phi-\psi\|,\qquad \forall t\geqslant0, \phi,\psi\in C(M). 
	\end{align*}
	In particular, if $L_u(x,v,u)\leqslant0$ for all $(x,v,u)\in TM\times\R$, then
	\begin{align*}
		\|T_t\phi-T_t\psi\|\leqslant \|\phi-\psi\|,\qquad \forall t\geqslant0, \phi,\psi\in C(M).
	\end{align*}
\end{enumerate}	
\end{Lem}

\section{Qualitative estimate of topological entropy}

\subsection{Lower bound}

Now we can estimate the topological entropy $h_{\rm top}(T_t)$. By the non-expansivity property of the semigroup $\{T_t\}$ in the case $L_u\leqslant0$, the following consequence can be easily obtained from Lemma \ref{lem:LO_semigroup} (4) and Lemma \ref{lem:entropy1}.

\begin{Pro}\label{pro:zero topology}
Under the assumption \emph{(L1)-(L3)}, if
\begin{align*}
	L_u(x,v,u)\leqslant0,\qquad\forall (x,v,u)\in TM\times\R,
\end{align*}
then $h_{\rm top}(T_t)=0$.
\end{Pro}

The following result provides a lower bound of $h_{\rm top}(T_t)$.

\begin{The}\label{thm:entrop2}
Under the assumptions \emph{(L1)-(L3)} and that there exists \(0<\lambda\leqslant \kappa\) such that
\begin{align*}
	\lambda\leqslant L_u(x,v,u)\leqslant \kappa , \qquad\forall (x,v,u)\in TM\times\R.
\end{align*}
\begin{enumerate}[\rm (1)]
	\item For any compact subset of the form $
		K=\Set{\phi=\phi_0+\delta}{\delta\in[0,a]}\subset C(M)$, where $\phi_0\in C(M)$ and $a>0$, it holds that $g(K)\geqslant\lambda$.
	\item $h_{\rm top}(T_t)\geqslant\lambda$.
\end{enumerate}
\end{The}

\begin{proof}
Fix $\varepsilon>0$ and $t>0$. Let $\phi_m=\phi_0+2m e^{-\lambda t}\varepsilon$, $m\in\N$. We claim
\begin{equation}\label{eq:dyn_ball1}
	\|T_t\phi_m-T_t\phi_{m-1}\|\geqslant 2\varepsilon,\qquad\forall m\in\N.
\end{equation}
Given $m\in\N$ and $x\in M$, let $\xi\in\mathcal{A}_{t,x}$ be a minimizer for $T_t\phi_m(x)$. That is 
\begin{align*}
	T_t\phi_m(x)=\phi_m(\xi(0))+\int^t_0L(\xi(s),\dot{\xi}(s),u_{\xi}(s))\ ds=u_{\xi}(t),
\end{align*}
where $u_{\xi}:[0,t]\to\R$ is determined by
\begin{align*}
	\begin{cases}
		\dot{u}_\xi=L(\xi,\dot{\xi},u_{\xi}),\\
		u_{\xi}(0)=\phi_k(\xi(0))=\phi_0(\xi(0))+2m e^{-\lambda t}\varepsilon.
	\end{cases}
\end{align*}
On the other hand, 
\begin{align*}
	T_t\phi_{m-1}(x)\leqslant \phi_{m-1}(\xi(0))+\int^t_0L(\xi(s),\dot{\xi}(s),v_{\xi}(s))\ ds=v_{\xi}(t), 
\end{align*}
with $v_{\xi}:[0,t]\to\R$ solving the equation
\begin{align*}
	\begin{cases}
		\dot{v}_\xi=L(\xi,\dot{\xi},v_{\xi}),\\
		v_{\xi}(0)=\phi_{m-1}(\xi(0))=\phi_0(\xi(0))+2(m-1)e^{-\lambda t}\varepsilon.
	\end{cases}
\end{align*}
Set $w_{\xi}=u_{\xi}-v_{\xi}$, then we have
\begin{align*}
	\begin{cases}
		\dot{w}_\xi=L(\xi,\dot{\xi},u_{\xi})-L(\xi,\dot{\xi},v_{\xi})=\widehat{L_u}\cdot w_{\xi},\\
		w_{\xi}(0)=u_{\xi}(0)-v_{\xi}(0)=2e^{-\lambda t}\varepsilon,
	\end{cases}
\end{align*}
where
\begin{align*}
	\widehat{L_u}(s)=\int^1_0L_u(\xi(s),\dot{\xi}(s),u_{\xi}(s)+\theta (u_{\xi}(s)-v_{\xi}(s)))\ d\theta,\quad s\in[0,t].
\end{align*}
It follows
\begin{align*}
	w_{\xi}(t)=e^{\int^t_0\widehat{L_u}\ ds}\cdot 2e^{-\lambda t}\varepsilon\geqslant e^{\lambda t}\cdot 2e^{-\lambda t}\varepsilon=2\varepsilon,
\end{align*}
and
\begin{align*}
	T_t\phi_{m-1}(x)\leqslant v_{\xi}(t)=u_{\xi}(t)-w_{\xi}(t)\leqslant T_t\phi_m (x)-2\varepsilon,
\end{align*}
which implies our claim \eqref{eq:dyn_ball1}.

Thus, for any $\phi\in C(M)$ the dynamical ball $B(\phi,t,\varepsilon)$ contains only a unique $\phi_m$ with $m\in\N\cup\{0\}$. Then, 
\begin{align*}
	g_t(\varepsilon,K)\geqslant\#\{\phi_m \in K\}\geqslant\left\lfloor\frac a{2e^{\lambda t}\varepsilon}\right\rfloor+1\geqslant\frac a{2e^{\lambda t}\varepsilon}.
\end{align*}
Hence, we obtain
\begin{align*}
	g(K)=\lim_{\varepsilon\to0^+}\limsup_{t\to\infty}\frac 1t\log g_t(\varepsilon,K)\geqslant\lim_{\varepsilon\to0^+}\limsup_{t\to\infty}\frac 1t\log \frac a{2e^{\lambda t}\varepsilon} =\lambda,
\end{align*}
and this completes the proof of (1).

Item (2) is immediate from the definition of $h_{\rm top}(T_t)$ together with (1). 
\end{proof}

\subsection{Upper bound}

Let $M$ be a compact and connected manifold without boundary. A function $l=l(x,v):TM\to\R$ is called a Tonelli Lagrangian if it is of class $C^2$ and satisfies the following conditions:
\begin{enumerate}[\rm (l1)]
	\item $l_{vv}(x,v)>0$, $\forall (x,v)\in TM$;
	\item There exist a super-linear continuous function $\theta_0:[0,+\infty)\to[0,+\infty)$ and a constant $c_0\geqslant0$ such that
	\begin{align*}
		l(x,v)\geqslant\theta_0(|v|_x)-c_0,\qquad\forall(x,v)\in TM.
	\end{align*}
\end{enumerate}

In this section, we will estimate the topological entropy $h_{\rm top}(T_t)$ in the case
\begin{align*}
	L(x,v,u)=l(x,v)+\lambda u,\quad (x,v,u)\in TM\times\R,
\end{align*}
where $l=l(x,v):TM\to\R$ is a Tonelli Lagrangian and $\lambda>0$ is a constant. It is well known that this type of contact system is equivalent to a classical non-autonomous system, which has some properties in the following Lemma \ref{lem:non-auto}.

\begin{Lem}\label{lem:non-auto}
	Suppose $M$ is a compact manifold, $l(x,v):TM\to\R$ is a Tonelli Lagrangian, and $\lambda>0$ is a constant. For any $t>0$, $y,x\in M$, let
	\begin{align*}
		A(t,y,x)=\inf_{\xi\in\Gamma^{t}_{y,x}}\int_{0}^{t}e^{-\lambda s}l(\xi(s),\dot{\xi}(s))ds.
	\end{align*}
	\begin{enumerate}[\rm (1)]
		\item For any $t>0$, $y,x\in M$, there exists a minimizer $\xi$ for $A(t,y,x)$. Moreover, there exists $K_1>0$ such that if $t\geqslant 1$, we have
		\begin{align*}
			|\dot{\xi}(s)|\leqslant K_1,\quad \forall s\in[0,t].
		\end{align*}
		\item There exists $K_0>0$ such that for all $t\geqslant1$, $y\in M$, the function $A(t,y,\cdot)$ is $e^{-\lambda t}K_0$-Lipschitz.
		\item If $M=\T^d$ and $l(x,v)=l(v)$ is independent of $x$, then for all $t\geqslant1$, $x\in\T^d$, the function $A(t,\cdot,x)$ is $e^{-\lambda t}K_0$-Lipschitz.
	\end{enumerate}
\end{Lem}

\begin{proof}
The result in (1) is well known. One can refer to \cite{Fathi_book} for the proof.

For any \( x \in M \) and minimizer \( \xi \) for \( A(t,y,x) \), the result (1) implies
\begin{align*}
	|p(t)| = |e^{-\lambda t} l_v(\xi(t), \dot{\xi}(t))| \leqslant e^{-\lambda t} K_0.
\end{align*}
Combining this with the relation (see \cite{Cannarsa_Sinestrari_book})
\begin{align*}
	D^+_x A(t,y,x) = \text{co} \Set{ p(t) = e^{-\lambda t} l_v(\xi(t), \dot{\xi}(t))}{\xi \text{ is a minimizer}},
\end{align*}
where $\text{``co''}$ stands for the convex hull, we conclude that \( A(t,y,\cdot) \) is \( e^{-\lambda t} K_0 \)-Lipschitz continuous.

Now, we turn to prove statement (3). Under this condition, for any \( y \in M \) and minimizer \( \xi \) for \( A(t,y,x) \), the Euler-Lagrange equation implies
\begin{align*}
	\dot{p}(s) = e^{-\lambda t} l_x(\dot{\xi}(s)) = 0, \quad \forall s \in [0,t].
\end{align*}
This leads to \( |p(0)| = |p(t)| \leqslant e^{-\lambda t} K_0 \). Additionally, we have
\begin{align*}
	D^+_y A(t,y,x) = - \text{co} \Set{ p(0) = e^{-\lambda t} l_v(\xi(0), \dot{\xi}(0)) }{\xi \text{ is a minimizer}}.
\end{align*}
 Therefore, we conclude that \( A(t, \cdot, x) \) is \( e^{-\lambda t} K_0 \)-Lipschitz continuous.
\end{proof}

\begin{Pro}\label{pro:dis lip}
Suppose $M$ is a compact manifold, $l(x,v):TM\to\R$ is a Tonelli Lagrangian, $\lambda>0$ is a constant. Let
\begin{align*}
	L(x,v,u):TM\times\R\to\R,\quad L(x,v,u)=l(x,v)+\lambda u.
\end{align*}
\begin{enumerate}[\rm (1)]
	\item For any $\varphi\in C^0(M)$, $t>0$, $x\in M$, we have
	\begin{align*}
		T_t\varphi(x)=\inf_{y\in M}\{e^{\lambda t}\varphi(y)+D(t,y,x)\},
	\end{align*}
	where
	\begin{align*}
		D(t,y,x)=\inf_{\xi\in\Gamma^{t}_{y,x}}\int_{0}^{t}e^{\lambda(t-s)}l(\xi(s),\dot{\xi}(s))ds.
	\end{align*}
	\item If $M=\T^d$ and $l(x,v)=l(v)$ is independent of $x$, then there exists $K_0>0$ such that for all $t\geqslant1$, $x\in\T^d$, the function $D(t,\cdot,x)$ is $K_0$-Lipschitz continuous.
\end{enumerate}
\end{Pro}

\begin{proof}
By definition, we have
\begin{align*}
	T_t \phi(x) = \inf_{y \in M} \inf_{\xi \in \Gamma_{y,x}^{t}} \left( \phi(y) + \int_0^t\big\{l(\xi(s), \dot{\xi}(s)) + \lambda u_{\xi}(s)\big\}\, ds \right)
\end{align*}
where \( u_{\xi}(s) \) satisfies the differential equation
\begin{equation}\label{eq:cara dis}
	\begin{cases}
		\dot{u}_{\xi}(s) = l(\xi(s), \dot{\xi}(s)) + \lambda u_{\xi}(s), & \quad s \in [0,t], \\
		u_{\xi}(0) = \phi(y). &
	\end{cases}
\end{equation}
Solving \eqref{eq:cara dis}, we obtain
\begin{align*}
	u_{\xi}(t) = e^{\lambda t} \phi(y) + \int_0^t e^{\lambda(t-s)} l(\xi(s), \dot{\xi}(s)) \, ds.
\end{align*}
Thus, we have
\begin{align*}
	T_t \phi(x) &= \inf_{y \in M} \inf_{\xi \in \Gamma_{y,x}^{t}} \left( e^{\lambda t} \phi(y) + \int_0^t e^{\lambda(t-s)} l(\xi(s), \dot{\xi}(s)) \, ds \right) \\
	&= \inf_{y \in M} \left( e^{\lambda t} \phi(y) + D(t,y,x) \right),
\end{align*}
where \( D(t,y,x) \) is the integral term. This completes the proof of (1).

To prove (2), for any \( t \geqslant 1 \), \( x, y \in \mathbb{T}^d \), we have
\begin{align*}
	D(t,y,x) = \inf_{\xi \in \Gamma^{t}_{y,x}} e^{\lambda t} \int_0^t e^{-\lambda s} l(\xi(s), \dot{\xi}(s)) \, ds = e^{\lambda t} A(t,y,x).
\end{align*}
Combining this with Lemma \ref{lem:non-auto} (3), we know that \( D(t,y,x) \) is \( K_0 \)-Lipschitz with respect to \( y \).
\end{proof}


For any compact subset $V\subset C^0(\T^d)$, we try to get a finite estimate of $g_t(\varepsilon,V)$. That is, we need to construct finite functions $\psi\in C^0(\T^d)$ such that
\begin{align*}
	\bigcup_{\psi}B(\psi,t,\varepsilon)\supset V.
\end{align*}
So we hope $T_t\psi$ is determined by the values of $\psi$ at finite points $z_i$'s, see the following Lemma \ref{lem:finite point}.

\begin{Lem}\label{lem:finite point}
	Under the assumption of Proposition \ref{pro:dis lip} (2), suppose $\{z_i\}_{i=1}^{m}$ is a set of finite points on $\T^d $. For any function $\psi:\{z_i\}_{i=1}^{m}\to\R$, let
	\begin{align*}
		&K_{\psi}=\max\{K_0,\mbox{\rm Lip}\,[\psi]\},\\
		&\tilde{\psi}:\T^d \to\R,\quad \tilde{\psi}(x)=\min_{i=1,\ldots,m}\psi(z_i)+K_{\psi}d(z_i,x).
	\end{align*}
	\begin{enumerate}[\rm (1)]
		\item $\tilde{\psi}$ is $K_{\psi}$-Lipschitz, and
		\begin{align*}
			\tilde{\psi}(z_i)=\psi(z_i),\quad i=1,\ldots,m.
		\end{align*}
		\item We have
		\begin{align*}
			T_t^-\tilde{\psi}(x)=\min_{i=1,\ldots,m}e^{\lambda t}\psi(z_i)+D(t,z_i,x),\quad \forall t\geqslant1,\ x\in\T^d .
		\end{align*}
	\end{enumerate}
\end{Lem}

\begin{proof}
Given \( x, y \in M \), there exists \( z_i \) such that
\begin{align*}
	\tilde{\psi}(x) &= \psi(z_i) + K_{\psi} \, d(z_i, x), \\
	\tilde{\psi}(y) &\leqslant \psi(z_i) + K_{\psi} \, d(z_i, y).
\end{align*}
This implies
\begin{align*}
	\tilde{\psi}(y) - \tilde{\psi}(x) &\leqslant K_{\psi} \, d(z_i, y) - K_{\psi} \, d(z_i, x) \leqslant K_{\psi} \, d(x, y).
\end{align*}
Similarly, we have \( \tilde{\psi}(x) - \tilde{\psi}(y) \leqslant K_{\psi} \, d(x, y) \), which shows that \( \tilde{\psi} \) is a Lipschitz continuous function with Lipschitz constant \( K_{\psi} \). Recalling that \( K_{\psi} \geqslant \mathrm{Lip}[\psi] \), for any \( i, j = 1, \ldots, m \), we have
\begin{align*}
	\psi(z_i) &\leqslant \psi(z_j) + \text{Lip}(\psi) \cdot d(z_j, z_i) \\
	&\leqslant \psi(z_j) + K_{\psi} \cdot d(z_j, z_i).
\end{align*}
It follows that \( \psi(z_i) = \tilde{\psi}(z_i) \) for all \( i = 1, \ldots, m \).

Now, we turn to the proof of (2). For any \( y \in \mathbb{T}^d \), there exists \( z_i \) such that \( \tilde{\psi}(y) = \psi(z_i) + K_{\psi} d(z_i, y) \). Now we have
\begin{align*}
	& e^{\lambda t} \psi(z_i) + D(t, z_i, x) \\
	\leqslant & \, e^{\lambda t} (\tilde{\psi}(y) - K_{\psi} d(z_i, y)) + D(t, y, x) + K_0 d(z_i, y) \\
	\leqslant & \, e^{\lambda t} \tilde{\psi}(y) + D(t, y, x), \quad \forall t \geqslant 1, \ x \in \mathbb{T}^d.
\end{align*}
This together with (1) implies
\begin{align*}
	T_t \tilde{\psi}(x) &= \inf_{y \in \mathbb{T}^n} e^{\lambda t} \psi(y) + D(t, y, x) \\
	&= \min_{i = 1, \ldots, m} e^{\lambda t} \psi(z_i) + D(t, z_i, x), \quad \forall t \geqslant 1, \ x \in \mathbb{T}^d.
\end{align*}
\end{proof}

\begin{The}\label{thm:dis}
	Under the assumption of Proposition \ref{pro:dis lip} (2), for any compact subset $V\subset C^0(\T^d )$, there holds $g(V)\leqslant\lambda$. Moreover, we have $h_{\rm top}(T_t)=\lambda$.
\end{The}

\begin{proof}
By compactness of $V$, there exists $R>0$ such that
\begin{align*}
	\|\phi\|_{C^0}\leqslant R,\quad \forall \phi\in V.
\end{align*}
Fix any $\varepsilon>0$, there exist finite points $\{z_i\}_{i=1}^m$ such that
\begin{align*}
	\bigcup_{i=1}^{m}B(z_i,\varepsilon)=\T^d .
\end{align*}
For any $n\in\N$, $k=0,1,\ldots,\lfloor\frac{2Re^{\lambda n}}{\varepsilon}\rfloor$, let
\begin{align*}
	U_{n,k}=\Set{-R+ke^{-\lambda n}\varepsilon+l\cdot e^{-\lambda j}\varepsilon}{j=1,\ldots,n,\ l=0,1,\ldots,\left\lceil\frac{K_0}{\varepsilon}\right\rceil+1}.
\end{align*}
Consider the set of functions
\begin{align*}
	&F_n=\Set{\psi:\{z_i\}_{i=1}^m\to U_{n,k}}{\ k=0,1,\ldots,\left\lfloor\frac{2Re^{\lambda n}}{\varepsilon}\right\rfloor},\\
	&\tilde{F}_n=\Set{\tilde{\psi}}{\psi\in F_n},
\end{align*}
where $\tilde{\psi}$ is defined by Lemma \ref{lem:finite point}. We claim that
\begin{equation}\label{eq:cover of V}
	\bigcup_{\tilde{\psi}\in\tilde{F}_n}\tilde{B}(\tilde{\psi},n,(1+K_0)\varepsilon)\supset V,
\end{equation}
where
\begin{align*}
	\tilde{B}(\tilde{\psi},n,(1+K_0)\varepsilon)=\Set{\phi\in C^0(\T^d )}{\|T_k\varphi-T_k\tilde{\psi}\|_{C^0}<(1+K_0)\varepsilon,\ \forall k=1,\ldots,n}.
\end{align*}
The relation \eqref{eq:cover of V} implies
\begin{align*}
	\tilde{g}_n(V,(1+K_0)\varepsilon)&\leqslant\#(\tilde{F}_n)=\#(F_n)=\sum_{k=0}^{\left\lfloor\frac{2Re^{\lambda n}}{\varepsilon}\right\rfloor}\#^m(U_{n,k})\\
	&\leqslant \left(\frac{2Re^{\lambda n}}{\varepsilon}+1\right)n^{m}\left(\frac{K_0}{\varepsilon}+3\right)^m,
\end{align*}
where
\begin{align*}
	\tilde{g}_n(V,(1+K_0)\varepsilon)=\min\Set{\#(E)}{E\subset C^0(\T^d ),\bigcup_{\psi\in E}\tilde{B}(\psi,n,(1+K_0)\varepsilon)\supset V}.
\end{align*}
It follows from Lemma \ref{lem:LO_semigroup} (4), Proposition \ref{pro:discrete} and Proposition \ref{pro:cut} that
\begin{align*}
	g(V)&=\lim_{\varepsilon\to 0^+}\limsup_{n\to+\infty}\frac{1}{n}\log \tilde{g}_n(V,\varepsilon)\\
	&=\lim_{\varepsilon\to 0^+}\limsup_{n\to+\infty}\frac{1}{n}\log \tilde{g}_n(V,(1+K_0)\varepsilon)\\
	&\leqslant \lim_{\varepsilon\to 0^+}\limsup_{n\to+\infty}\frac{1}{n}\log \left(\frac{2Re^{\lambda n}}{\varepsilon}+1\right)n^{m}\left(\frac{K_0}{\varepsilon}+3\right)^m\\
	&=\lambda.
\end{align*}
By the arbitrariness of \( V \), we have
\begin{align*}
	h_{\rm top}(T_t) = \sup_{V \subset C^0(\mathbb{T}^n)} g(V) \leqslant \lambda.
\end{align*}
Combining this with Theorem \ref{thm:entrop2}, we conclude that 
\begin{align*}
	h_{\rm top}(T_t) = \lambda.
\end{align*}

Finally, we only need to prove the claim \eqref{eq:cover of V}. In fact, for any $\phi\in V$, there exists a unique $k\in\{0,1,\ldots,\lfloor\frac{2Re^{\lambda n}}{\varepsilon}\rfloor\}$ such that
\begin{align*}
	-R+ke^{-\lambda n}\varepsilon\leqslant \min_{x\in\T^d }\phi(x) < -R+(k+1)e^{-\lambda n}\varepsilon.
\end{align*}
Let $\psi:\{z_i\}_{i=1}^{m}\to\R$,
\begin{align*}
	\psi(z_i)=\max\Set{c\in U_{n,k}}{c\leqslant\phi(x),\forall x\in B(z_i,\varepsilon)}.
\end{align*}
Then $\psi\in F_n$ and $\tilde{\psi}\in\tilde{F}_n$. Now we only need to prove
\begin{align*}
	\phi\in \tilde{B}(\tilde{\psi},n,(1+K_0)\varepsilon),
\end{align*}
that is,
\begin{equation}\label{eq:cover phi}
	\|T_j\tilde{\psi}-T_j\phi\|_{C^0}<(1+K_0)\varepsilon,\quad \forall j=1,\ldots,n.
\end{equation}
For any $j\in\{1,\ldots,n\}$, $x\in\T^d $, there exists $y_x\in\T^d $ such that
\begin{align*}
	T_j\phi(x)=e^{\lambda j}\phi(y_x)+D(j,y_x,x).
\end{align*}
Choose $z_i$ such that $y_x\in B(z_i,\varepsilon)$. By the choice of $\psi$, we have $\psi(z_i)\leqslant\phi(y_x)$. Using Proposition \ref{pro:dis lip}, we know that
\begin{equation}\label{eq:pf 1}
\begin{split}
T_j\tilde{\psi}(x)&\leqslant e^{\lambda j}\psi(z_i)+D(j,z_i,x)\\
&\leqslant e^{\lambda j}\phi(y_x)+D(j,y_x,x)+(D(j,z_i,x)-D(j,y_x,x))\\
&\leqslant T_j\phi(x)+K_0\varepsilon.
\end{split}
\end{equation}
On the other hand, the choice of $k$ and $\psi$ implies that there exists $z_{i_0}$ such that
\begin{align*}
	\psi(z_{i_0})=\min_{1=1,\ldots,m}\psi(z_i)=-R+ke^{-\lambda n}\varepsilon.
\end{align*}
Lemma \ref{lem:finite point} implies that there exists $z_{i_x}$ such that
\begin{align*}
	T_j\tilde{\psi}(x)=e^{\lambda j}\psi(z_{i_x})+D(j,z_{i_x},x)\leqslant e^{\lambda j}\psi(z_{i_0})+D(j,z_{i_0},x),
\end{align*}
which leads to
\begin{equation}\label{eq:choiceofy}
\begin{aligned}
	\psi(z_{i_x})&\leqslant\psi(z_{i_0})+e^{-\lambda j}(D(j,z_{i_0},x)-D(j,z_{i_x},x))\\
	&\leqslant -R+ke^{-\lambda n}\varepsilon+K_0 e^{-\lambda j}.
\end{aligned}
\end{equation}
Now, by the definition of \( U_{n,k} \) and the choice of \( \psi \), we know that there exists \( y \in B(z_{i_x}, \varepsilon) \) such that
\begin{align*}
	\psi(z_{i_x}) \leqslant \phi(y) < \psi(z_{i_x}) + e^{-\lambda j} \varepsilon.
\end{align*}
 If it were not the case, Pigeonhole principle yields \[ \psi(z_{i_x})\geq -R + k e^{-\lambda n}\varepsilon  + \left( \left\lceil \frac{K_0}{\varepsilon}\right\rceil +1\right) \cdot e^{-\lambda j}\varepsilon>-R + k e^{-\lambda n}\varepsilon  + K_0 e^{-\lambda j},\] which contradicts (\ref{eq:choiceofy}). 
Thus we have that
\begin{equation}\label{eq:pf 2}
	\begin{split}
		T_j\tilde{\psi}(x)&=e^{\lambda j}\psi(z_{i_x})+D(j,z_{i_x},x)\\
		&\geqslant e^{\lambda j}(\phi(y)-e^{-\lambda j}\varepsilon)+D(j,y,x)+(D(j,z_{i_x},x)-D(j,y,x))\\
		&\geqslant e^{\lambda j}\phi(y)-\varepsilon+D(j,y,x)-K_0 d(z_{i_x},y)\\
		&\geqslant T_j\phi(x)-(1+K_0)\varepsilon.
	\end{split}
\end{equation}
It is easy to see that \eqref{eq:cover phi} follows directly from \eqref{eq:pf 1} and \eqref{eq:pf 2}. This completes the proof of our claim \eqref{eq:cover of V}.
\end{proof}

\bibliographystyle{plain}
\bibliography{mybib.bib}
\end{document}